\documentclass[microtype]{gtpart}

\usepackage{latexsym,enumitem}
\usepackage{amssymb}
\usepackage[cp850]{inputenc}
\usepackage{epsfig}
\usepackage{psfrag}
\usepackage{amsthm}
\usepackage{amscd}
\usepackage{amsmath}
\usepackage{amsfonts}
\usepackage{graphics,caption}
\usepackage[all]{xy}
\usepackage{etoolbox}
\patchcmd{\quote}{\rightmargin}{\leftmargin 2em \rightmargin}{}{}
\usepackage[sort,numbers]{natbib} 
\usepackage{mathtools}
\usepackage{enumitem}

\usepackage{overpic}

\usepackage[english]{babel}
\usepackage{comment}
\usepackage{color}
\usepackage{hyperref}

\let\phi\varphi

\let\epsilon\varepsilon
\let\subset\subseteq
\newcommand{\be}{\begin{equation*}}
 \newcommand{\ee}{\end{equation*}}
\newcommand{\bpf}{\begin{dimo}}
\newcommand{\epf}{\end{dimo}}
\newcommand{\bdefi}{\begin{defin}}
\newcommand{\edefi}{\end{defin}}
\newcommand{\bthm}{\begin{thm}}
\newcommand{\ethm}{\end{thm}}
\newcommand{\blem}{\begin{lem}}
\newcommand{\elem}{\end{lem}}
\newcommand{\bcor}{\begin{cor}}
\newcommand{\ecor}{\end{cor}}

\newcommand{\bprop}{\begin{prop}}
\newcommand{\eprop}{\end{prop}}
\newcommand{\bese}{\begin{ese}}
\newcommand{\eese}{\end{ese}}
\newcommand{\brem}{\begin{rem}}
\newcommand{\erem}{\end{rem}}
\newcommand{\bpfc}{\begin{dimoclaim}}
\newcommand{\epfc}{\end{dimoclaim}}
\newcommand{\eps}{\epsilon}
\newcommand{\al}{\alpha}
\newcommand{\sig}{\sigma}

\newcommand{\rar}{\rightarrow} 

\newcommand{\abs}[1]{\left\lvert#1\right\rvert}						
\newcommand{\set}[1]{\left\{#1\right\}}					
	
\newcommand{\quotient}[2]{\left.\raisebox{.1em}{$#1\!$}\middle/\raisebox{-.1em}{$#2$}\right.}

\DeclareMathOperator{\bR}{\mathbb R}			
\DeclareMathOperator{\bZ}{\mathbb Z}			

\DeclareMathOperator{\im}{Im} 

\newenvironment{quot}
{
	\vspace{-0.2cm}
	\vspace{0.2cm}
}

\theoremstyle{definition}
\newtheorem{d1}{Definition}

\newenvironment{defin}
{
	\begin{quot}
		\begin{d1}
		}
		{\end{d1}
	\end{quot}

}

\theoremstyle{definition}
\newtheorem{r1}[d1]{Remark}

\newenvironment{rem}
{
	\begin{quot}
		\begin{r1}
		}
		{\end{r1}
	\end{quot}
}

\theoremstyle{definition}
\newtheorem{e1}[d1]{Exercise}

\theoremstyle{definition}
\newtheorem{ese1}[d1]{Example}

\newenvironment{ese}
{
	\begin{quot}
		\begin{ese1}
	}
	{	
		\end{ese1}
	\end{quot}
}

\theoremstyle{definition}

\theoremstyle{definition}
\newtheorem{f2}[d1]{Fact}

\theoremstyle{definition}

\theoremstyle{definition}

\theoremstyle{definition}
\newtheorem{t1}[d1]{Theorem}

\newenvironment{thm}
{
	\begin{quot}
		\begin{t1}}
		{\end{t1}
	\end{quot}
}

\theoremstyle{definition}
\newtheorem*{T1*}{Theorem}

\newenvironment{teor*}
{
	\begin{quot}
		\begin{T1*}}
		{\end{T1*}
	\end{quot}
}

\newenvironment{dimo}
{\begin{proof}[Proof]
	}
	{\end{proof}}

\newenvironment{dimoclaim}{\emph{Proof of Claim:}\;}{\hfill$\square$}

	\theoremstyle{definition}
	\newtheorem{l1}[d1]{Lemma}
	
	\newenvironment{lem}
	{
		\begin{quot}
			\begin{l1}}
			{\end{l1}
		\end{quot}
	}
	\theoremstyle{definition}
	\newtheorem{p1}[d1]{Proposition}
	
	\newenvironment{prop}
	{
		\begin{quot}
			\begin{p1}}
			{\end{p1}
		\end{quot}
	}
	
	\theoremstyle{definition}
	\newtheorem{c1}[d1]{Corollary}
	
	\newenvironment{cor}
	{
		\begin{quot}
			\begin{c1}}
			{\end{c1}
		\end{quot}
	}

		\renewenvironment{abstract}
	{\list{}{\rightmargin\leftmargin}%
		\item[\textbf{Abstract:}]\relax}
	{\endlist}

 \newtheorem*{Theorem*}{Theorem}
 \newtheorem*{Proposition*}{Proposition}
 \newtheorem*{Lemma*}{Lemma}

\usepackage{amsmath,amssymb,amsthm}
\usepackage{thmtools}
\usepackage{thm-restate}

\makeatletter
\DeclareFontFamily{OMX}{MnSymbolE}{}
\DeclareSymbolFont{MnLargeSymbols}{OMX}{MnSymbolE}{m}{n}
\SetSymbolFont{MnLargeSymbols}{bold}{OMX}{MnSymbolE}{b}{n}
\DeclareFontShape{OMX}{MnSymbolE}{m}{n}{
    <-6>  MnSymbolE5
   <6-7>  MnSymbolE6
   <7-8>  MnSymbolE7
   <8-9>  MnSymbolE8
   <9-10> MnSymbolE9
  <10-12> MnSymbolE10
  <12->   MnSymbolE12
}{}
\DeclareFontShape{OMX}{MnSymbolE}{b}{n}{
    <-6>  MnSymbolE-Bold5
   <6-7>  MnSymbolE-Bold6a
   <7-8>  MnSymbolE-Bold7
   <8-9>  MnSymbolE-Bold8
   <9-10> MnSymbolE-Bold9
  <10-12> MnSymbolE-Bold10
  <12->   MnSymbolE-Bold12
}{}

\let\llangle\@undefined
\let\rrangle\@undefined
\DeclareMathDelimiter{\llangle}{\mathopen}%
                     {MnLargeSymbols}{'164}{MnLargeSymbols}{'164}
\DeclareMathDelimiter{\rrangle}{\mathclose}%
                     {MnLargeSymbols}{'171}{MnLargeSymbols}{'171}
\makeatother

\newcommand{\Diffeo}{\operatorname{Diffeo}}
\newcommand{\Mod}{\operatorname{Mod}}
\newcommand{\bS}{\mathbb{S}}
\newcommand{\bD}{\mathbb{D}}

\renewcommand{\hat}{\widehat}
\newcommand{\wt}{\widetilde}

\renewcommand{\setminus}{\smallsetminus}

\begin{document}
\title[Knots in circle bundles are determined by their complements ]{Knots in circle bundles are determined by their complements}
\author{Tommaso Cremaschi}
\email{tommaso.cremaschi@uni.lu}

\author{Andrew Yarmola}
\email{andrew.yarmola@yale.edu}
\thanks{The first author acknowledges support from XXX.}

\date{v1, \today}

\maketitle
\begin{abstract}
We resolve a case of the oriented knot complement conjecture by showing that knots in an orientable circle bundle $N$ over a genus $g \geq 2$ surface $S$ are determined by their complements. We apply this to the setting of canonical knots in the unit and projective tangent bundles, which are knots that are the set of tangents to a closed curve on $S$. We show that canonical knots have homeomorphic complements if and only if their shadows differ by Reidemeister moves, (de)stabilizations, loops/cusps added by transvections, and mapping classes of $S$.
\end{abstract}

\section{Introduction}
The results of this paper are motivated by the {\it oriented knot complement conjecture} of \cite[1.81(D)]{kirby1997problems} and \cite[Conjecture 6.2]{gordon1990dehn}. The conjecture states that if two knots in an oriented, closed 3-manifold $N$ have orientation-homeomorphic exteriors, which are not solid tori, then there is orientation-preserving homeomorphism of $N$ taking one knot to the other. This conjecture was answered in the affirmative for knots in the 3-sphere by Gordon and Luecke \cite{GL1989} and for knots in $\mathbb{S}^2 \times \mathbb{S}^1$ and torus bundles over $\mathbb{S}^1$ by Gabai in \cite{Gabai}. Our main result confirms this conjecture for $N$ an orientable circle bundle over an orientable surface of negative Euler characteristic. In particular, we show that the homeomorphism of the exteriors extends to that of $N$ for non-trivial knots.

\begin{restatable}{theorem}{main}\label{SFone} Let $M$ be a compact orientable 3-manifold where $\partial M = T$ is a torus. Consider slopes $s_1$ and $s_2$ on $T$. Assume that $M(s_1) \cong^+ M(s_2)$ and that $N = M(s_1)$ is an orientable circle bundle over an orientable surface of negative Euler characteristic. Then either $s_1 = s_2$ or both of the induced embeddings of $M$ in $N$ are trivial.
\end{restatable}

We say that an embedding of manifolds $\psi:M\hookrightarrow N$ is \emph{trivial} if the image $\psi(M)$ is the exterior of an unknot in $N$. As an immediate corollary, we have:

\bcor Let $N$ be an orientable circle bundle over an orientable surface of negative Euler characteristic. Two knots $K_1, K_2$ in $N$ have homeomorphic exteriors if and only if $h(K_1) = K_2$ for some homeomorphism $h$ of $N$, possibly orientation reversing.
\ecor

Analogues of the above statements, often with exceptions and restrictions on the knots, are also known in several settings. A brief overview includes:

\begin{itemize}
\item knots in solid tori \cite{BERGE19911, GABAI19891, GABAI1990221}, with 0-bridge and 1-bridge knots as exceptions.
\item knots with Seifert fibered exteriors in Seifert fibered spaces \cite{rong1993some}.
\item non-hyperbolic knots in lens spaces and closed, atoroidal and irreducible  Seifert-fibered spaces \cite{matignon2010knot, ito_2022}.
\item null-homologous knots knots in some $L$-spaces and rational homology spheres \cite{gainullin2018heegaard, ito_2022}.
\item hyperbolic, homologically non-trivial knots in non-hyperbolic homology lens spaces \cite{Ichihara-Saito}.
\end{itemize}

Most recent results listed above make use of Heegaard Floer machinery. In contrast, our techniques rely heavily on topology and group-theoretic arguments, see Section \ref{secGL} for an outline of the proof.

{\bf Applications to knots in tangent bundles.} In recent years there has been a lot of interest in knots and links in tangent bundles. A particularly rich examples are realized as periodic orbits of some geodesic flow on the unit tangent bundle of a surface, see  \cite{CRMY2022,RM2020,CKMP2021,CRM2020}. Classically, Ghys \cite{G2007} showed that the periodic orbits of the geodesic flow over the modular surface $\Sigma_\text{mod}$ correspond to Lorenz links in $\mathbb S^3$, which are periodic orbits of the Lorenz attractor. Lorenz links enjoy lots of nice properties, they are prime, positive, fibered, hence amphicherical, and moreover both the genus and their braid index can be studied combinatorially, see \cite{BW1983} for more details.

 A motivating application of Theorem \ref{SFone} is to show such orbits are determined by their complements, extending results of \cite{RM2020}. Let $S$ be a surface and fix a collection $\kappa$ of closed (oriented) curves on $S$. Consider the set of tangents $\hat\kappa$ to $\kappa$ in the unit tangent bundle $UT(S)$ or the projective tangent bundle $PT(S)$. We call these \emph{canonical links}, or \emph{knots} when $\kappa$ is just one curve. We require that $\kappa$ has no self-tangencies, so that the lift is embedded. We will often refer to $\kappa$ as a \emph{diagram} or \emph{shadow} for $\hat\kappa$. In the case of $UT(S)$, we also assume that $\kappa$ is oriented and smooth.
For $PT(S)$, we will allow $\kappa$ to be piecewise smooth where at the non-smooth points of $\kappa$ must have cusps, i.e. have a uniquely defined tangent line. We call such links $\kappa$ \emph{cusp-smooth}. Note, smooth curves are also considered to be cusp-smooth. We use $M \smallsetminus L$ to denote the complement of an open regular neighborhood of a link $L$ in a 3-manifold $M$. We call $UT(S) \smallsetminus \hat\kappa$ and $PT(S) \smallsetminus \hat\kappa$ {\it canonical link exteriors}.

\brem If one works in the projective cotangent bundle and lifts to functionals whose kernel is perpendicular to a cusp-smooth curve, this gives a link tangent to the canonical contact structure on the projective cotangent bundle. Such links are called Legendrian. Topologically, these are the same as lifting tangents.
\erem

Our first observation is that any link in $UT(S)$ or $PT(S)$ is isotopic to a canonical one. 

\begin{restatable}{theorem}{links}\label{A} Let $S$ be a surface. For any link $L \subset UT(S)$ or $PT(S)$, there is a collection $\kappa$ of (cusp-)smooth closed (oriented) curves on $S$ such that $L$ is ambiently isotopic to $\hat\kappa$.
\end{restatable}

 Applying the above proposition to links in $\bD^2$, we see that 

\begin{restatable}{corollary}{linkscor}\label{CA} Every closed orientable connected 3-manifold $M$ is an $\mathbb{Z} \cup \{1/0\}$-Dehn filling of a canonical link exterior.\end{restatable}

We will consider two closed curves on a surface equivalent if they are isotopic up to flipping orientation. We say that two curves are equivalent up to some operations whenever a composition of such operations produces equivalent curves. In particular, since we are dealing with knots in manifolds with non-trivial mapping class groups, we will need to consider the action of this group.

\begin{restatable}{theorem}{appl}\label{GL} Let $S_1, S_2$ be orientable closed surfaces with $\chi(S_1) < 0$. Consider closed curves $\kappa_1 \subset S_1$ and $\kappa_2 \subset S_2$. Then, we have the following cases. \begin{enumerate}
\item[(a)] {\bf Case $\kappa_i$ smooth.} $UT(S_1) \smallsetminus \hat\kappa_1 \cong^+ UT(S_2) \smallsetminus \hat\kappa_2$ if and only if $S = S_1 \cong S_2$ and $\kappa_1, \kappa_2$ are equivalent up to $\Diffeo^+(S)$, loop Reidemeister moves and (de)stabilizations, and loops added by transvections.
\item[(b)] {\bf Case $\kappa_i$ cusp-smooth.} $PT(S_1) \smallsetminus \hat\kappa_1 \cong^+ PT(S_2) \smallsetminus \hat\kappa_2$ if and only if $S = S_1 \cong S_2$ and $\kappa_1, \kappa_2$ are equivalent up to $\Diffeo^+(S)$, Legendrian Reidemeister and (de)stabilizations, and cusps added by transvections.
\end{enumerate}
\end{restatable}

Theorem \ref{GL} generalizes the work of \cite{RM2020}, where it is proven for $UT(S_i) \smallsetminus\hat \kappa_i$ with the assumptions that $S_1 = S_2$, the $\kappa_i$ are smooth and minimal position, and $[\kappa_i]\neq 0$ in $H_1(S_1,\R)$ for $i = 1, 2$.
We will define Legendrian and loop Reidemeister moves, (de)stabilizations and transvections in Section \ref{reiddiffcan}. The later are fibre preserving mapping classes of the ambient circle bundle that twist the knot around the fibre. This corresponds to adding loops or cusps to the diagrams on $S$ in a controlled way. 

\paragraph{Organization.} In Section \ref{bkg}, we provide some background on circle bundles and discuss conventions. In Section  \ref{secGL}, we prove Theorem \ref{SFone} by first proving a similar result, Theorem \ref{SF}, for manifolds with multiple torus boundary components. Finally, in Section \ref{appsec} we study canonical links. 

\paragraph{Acknowledgements.} The first author acknowledges support from MSCA 101107744-DefHyp. We thank Jos{\'e}~Andr{\'e}s Rodr{\'\i}guez-Migueles for many useful conversations.

\section{Background}\label{bkg}
A properly embedded surface $S$ in a 3-manifold $M$ is \emph{incompressible} if it is $\pi_1$-injective and is said to be $\partial$-\emph{incompressible} if there is no disk $D\subset \overline {M\setminus S}$ such that $\partial D=\alpha\cup\beta$ with $\alpha\subset S$ and $\beta\subset\partial M$.

A connected, incompressible, $\partial$-incompressible, and non-$\partial$-parallel properly embedded surface in $3$-manifold will be called \emph{essential}. 

Given a class $c\in H_2(M)$ or $H_2(M, U \subset \partial M)$, its \emph{Thurston norm} is $x(c)=\min_{S} \chi_-(S)$ where $S$ is a properly embedded representative of $c$ and $\chi_-(S)=\max\set{-\chi(S), 0}$.

\subsection*{ Mapping class group of $UT(S)$ and $PT(S)$}

For an oriented manifold $N$, let $\Diffeo^+(N)$ denote the group of orientation preserving diffeomorphisms of $N$ and let $\Diffeo_0^+(N)$ denote the connected component of the identity. The mapping class group of $N$ is $\Mod(N) = \pi_0(\Diffeo^+(N)) = \Diffeo^+(N)/\Diffeo_0^+(N)$. Assume that $N$ is an oriented circle bundle over an orientable surface $S$ where $\chi(S) < 0$. Then there is a short exact sequence
\[ 0 \to H_1(S, \partial S; \mathbb Z) \xrightarrow{\mathfrak t} \Mod(N) \xrightarrow{\mathfrak h} \Mod(S) \to 0.\]
Since $\chi(S) < 0$, it is well known that any diffeomorphism is isotopic to one that maps fibers to fibers, so $\mathfrak h$ arises as the map induced on the space of fibers, which is $S$. The map $\mathfrak t$ corresponds to diffeomorphisms that map each fibre to itself. The authors of \cite{BC2023} call such diffeomorphism \emph{transvections}.  Transvections correspond to picking a representative $\al = \sum a_i \gamma_i$ of $[\al] \in H_1(S, \partial S; \bZ)$ as a weighted collection of disjoint simple closed curves and simple proper arcs and doing $a_i$ twists along the torus or annulus above each $\gamma_i$. In Figure \ref{tbtwist}, you can see an example of the action of a transvection on a piece of a knot. When $N = UT(S)$ or $PT(S)$, the short exact sequence above splits via the natural action of a surface mapping class on the tangent bundle. In particular, any mapping class of $N$ in this setting is a transvection followed by the lift of a surface mapping class. In the Section \ref{appsec}, we look at how isotopy and transvections act on canonical link diagrams.

\section{Knots in circle bundles and their complements}\label{secGL} 

To show that knots in circle bundles over closed hyperbolizable surfaces are determined by their complements (up to the action of diffeomorphisms), we prove the following analogue of the Gordon-Luecke Theorem \cite{GL1989}. The proof will involve facts about surface groups and topological arguments. A key element will be that our circle bundles have unique, up to isotopy, fibrations, see \cite{Ha2007}.

Recall from the introduction that an embedding of manifolds $\psi:M\hookrightarrow N$ is \emph{trivial} if the image $\psi(M)$ is the exterior of an unknot in $N$. The goal of this section is to prove:

\main*

\textbf{Proof Outline.}\label{outline} To prove Theorem \ref{SFone}, we first prove a version in which the circle bundles have non-empty boundary, see Theorem \ref{SF}. This with-boundary-version of the result requires a technical condition where the two induced embeddings of $M$ must agree on the un-filled boundary.

The proof of Theorem \ref{SF}, is very topological and is divided by various cases that use the Gordon-Luecke Theorem to deal with non-trivial knots in 3-balls and extensively use results of Gabai \cite{Gabai} on the Thurston norm of surfaces in 3-manifolds and their fillings. 

A case of Theorem \ref{SFone} will follow from Theorem \ref{SF}. To do so, we show that there are vertical fibers $M(s_1)$ and $M(s_2)$ whose pull backs are homotopic, see Proposition \ref{fibre}. Then, by drilling it out we fall into the setting of Theorem \ref{SF} and complete the proof.

\brem One aspect of our proof approach is that it is possible to trace stronger variants of Theorem \ref{SFone} in certain cases. In particular, for many classes of knots, the knot exterior $M$ has the property that there is a unique slope whose filling is an oriented circle bundle over an orientable surface of negative Euler characteristic.
\erem

We begin our proof of Theorems \ref{SF} and \ref{SFone} with a quick, purely topological, Lemma.

\blem\label{type} Let $M$ be a compact orientable 3-manifold where $\partial M$ is a union of tori $T_0,\dotsc, T_k$. Let $s_1$ and $s_2$ be slopes on the first boundary component $T = T_0$ of $M$. If $M(s_1)$ and $M(s_2)$ are both orientable circle bundles over orientable surfaces, then the base surfaces are diffeomorphic.
\elem

\bpf Let $N_i = M(s_i)$ be a circle bundles over the orientable surface $S_i$ for $i = 1, 2$. Let $T=T_0$ be the first boundary component of $M$. Place a basepoint on $T$ and realize $s_i$ as $\sigma_i \in \pi_1(M)$ on $T$. Looking at homology, by Mayer-Vietoris, we see that 

\[H_1(M; \mathbb Z)/\langle [\sigma_i] \rangle \cong H_1(N_i; \mathbb Z).\] 

When $|\partial M| > 1$, we know that $H_1(N_i; \mathbb Z) = \mathbb{Z}^{2g_i + k_i-1} \oplus \bZ $, where $g_i$ and $k_i> 0$ are the genus and number of boundary components of $S_i$. When $|\partial M| = 1$,  we have $H_1(N_i; \mathbb Z) = \mathbb{Z}^{2g_i} \oplus \mathbb Z/e_i \mathbb{Z}$, where  $e_i$ is the Euler number of the bundle $N_i$. 

Assume $S_1$ has boundary. By construction, so does $S_2$ and the number of components is the same, i.e. $k_1 = k_2 > 0$. In this setting, $H_1(N_i; \mathbb Z)$ has no torsion, so $\langle [\sig_i] \rangle$ either kills the torsion of $H_1(M; \mathbb Z)$ or $H_1(M; \mathbb Z)$ is torsion-free and $\langle [\sig_i] \rangle$ is trivial or kills a $\mathbb Z$ summand. In all cases, this gives $|2g_1 - 2g_2| \leq 1$, which means $g_1 = g_2$. In the closed case, $\langle [\sig_i] \rangle$ can have the additional behavior of creating the Euler number torsion. However, all the cases still give $|2g_1 - 2g_2| \leq 1 $, so $g_1 = g_2$ as before.
\epf
\brem Lemma \ref{type} also holds when one considers circle bundles over non-orientable surfaces. Namely, if both surfaces are non-orientable the same homology computations work. Additionally, no circle bundle over a non-orientable surface can be homeomorphic to a circle bundle over an orientable surface since the fundamental group of the latter never has two-torsion.
\erem

The proof of Theorem \ref{SFone} relies on a technical reduction to the case where $N$ has multiple torus boundary components and certain properties of surface groups. This first part is Theorem \ref{SF} below, which we now prove.

\bthm\label{SF} Let $M$ be a compact orientable 3-manifold where $\partial M$ is a union of tori and $|\partial M| \geq 2$.
Let $s_1$ and $s_2$ be slopes on $T = T_0$ where $\partial M = T_0 \cup \cdots \cup T_k$ and consider the natural inclusions $\psi_i: M \to M(s_i)$ for $i = 1,2$. Assume that
\begin{enumerate}
\item $M(s_1)$ is an orientable circle bundle over an orientable surface $S$ with $\chi(S) < 0$;
\item there is an orientation preserving diffeomorphism $g: M(s_1) \to M(s_2)$ such that $\left.(\psi_2^{-1} \circ g \circ \psi_1)\right|_{T_1 \cup \cdots \cup T_k}$ is isotopic to the identity on $T_1 \cup \cdots \cup T_k$.
\end{enumerate}
Then, either the slopes are equal $s_1 = s_2$ or both $\psi_1$ and $\psi_2$ are trivial. \ethm

\begin{proof} We may assume that one of the $\psi_i$'s is non-trivial, since otherwise we are done.

To begin, notice that $|\partial M| \geq 2$ and the diffeomorphism $g$, imply that $M(s_1)$ and $M(s_2)$ are both the trivial circle bundle $N = S \times \bS^1$ where $S$ has non-empty boundary. Fix embeddings $\psi_i:M \to N$ arising from the Dehn fillings that realize $M$ as the complement of an open tubular neighborhood of the cores $\eta_i$ of the filling, for $i = 1,2$. By replacing $\psi_1$ with $g \circ \psi_1$ followed by an isotopy, we can further assume that $\left.(\psi_2^{-1} \circ \psi_1)\right|_{T_1 \cup \cdots \cup T_k}$ is the identity.

Let $T= T_0$ be the first boundary component of $M$. Place a basepoint  $x$ on $T$ and realize the slopes $s_i$ as $\mu_i \in \pi_1(M,x)$ on $T$, picking an orientation for each. Under $\psi_i$, $\mu_i$ becomes a meridian of $\eta_i$ and induces a positive orientation on $\eta_i$. Isotope $\eta_i$ so it is never tangent to the fibers of $N$ and push it, transverse to the fibres, to get a longitude $\ell_i$ on $T$ for each $i = 1, 2$. 

If $M$ is reducible, then since $N$ is irreducible, it follows that the reducing sphere cuts off a ball containing $T$. In particular $M$ is a connect sum of $N$ and a knot complement in $\bS^3$. This sphere must be unique up to isotopy and we can assume that the knot is not an unknot by non-triviality of at least one of the embeddings. Thus, by uniqueness of prime decompositions and the fact that we are dealing with a non-trivial knot in a 3-ball,  
 the Gordon-Luecke Theorem implies that $s_1 = s_2$. We may now assume $M$ is irreducible. That is, neither $\eta_1$ or $\eta_2$ are contained in a 3-ball in $N$.

We have that $\pi_1(N) = \pi_1(S) \times \mathbb{Z}$, so $\eta_i$ is homotopic to $\alpha_i t^{m_i}$ where $\alpha_i \in \pi_1(S)$ and $t$ generates the fibre. Note, we will abuse notation using $\al_i$ for both the homotopy class in $\pi_1(S)$ and the actual shadow of $\eta_i$ under the canonical projection of $N$ to $S$. We now deal with various cases depending on conditions on $m_i$ and $\alpha_i$.  After the first case below, we may also assume that $M$ is $\partial$-irreducible since $N$ is $\partial$-irreducible and any compression disk with boundary on $T$ would either make $\eta_i$ null-homologous (in particular $m_i = 0$) or, if that disk has meridional boundary, it produces a non-separating sphere in $N$, which does not exist.

The proof will now continue by cases analysing the structure of $\eta_i$ and its homotopy class $\alpha_i t^{m_i}$ in $\pi_1(N)$.

\textbf{(1) Case $m_1 = 0$ or $m_2 = 0$.} 

We can assume $m_1 = 0$. Since $N = S \times \mathbb{S}^1$ is trivial and $m_1 = 0$, there is a horizontal section $F$ of $S$ in $N$ that intersects $\eta_1$ zero times algebraically. Let $R = \psi_1^{-1}(F)$ in $M$. Using $\partial N$ to denote $\partial M \setminus T$, we have that $[R] \in H_2(M, \partial N; \bR)$ is non-zero. Let $x: H_2(M, \partial N; \bR) \to  \bR$ be the Thurston norm. We claim that $x([R]) \geq |\chi(S)| > 0$. Indeed, let $R'$ be a properly embedded surface in $M$ realizing the Thurston norm of $[R]$. A priori, we do not know if $R'$ and $T$ are disjoint, but we can edit $R'$ so this is true as follows. By construction $[\partial R'] \cap [T] = 0$ in homology, so we can replace pairs of oppositely oriented components of $R' \cap T$ by attaching tubes that run parallel to $T$ in a regular neighborhood of a fixed annulus containing $\partial R'$ on $T$, see Figure \ref{tubingpic}.  Note, we can guarantee embeddedness by attaching tubes inductively on neighboring pairs of components with opposite orientation. This gives a possibly new Thurston-norm minimizing surface $R''$ with $[R''] = [R'] = [R]$ in $H_2(M, \partial N; \bR)$. After Dehn filling along $s_1$ on $T$, we see that $|\chi(R'')| \geq x([F]) = |\chi(S)| > 0$ since $[\psi_1(R'')] = [F]$. Further, since incompressible surfaces in $N$ with negative Euler characteristic are horizontal, we know that equality $|\chi(R'')| = |\chi(S)|$ holds if and only if $\psi_i(R'')$ is isotopic to a section disjoint from $\eta_i$ for both $i = 1,2$.

\begin{figure}
\begin{center}
\begin{overpic}[scale=0.9]{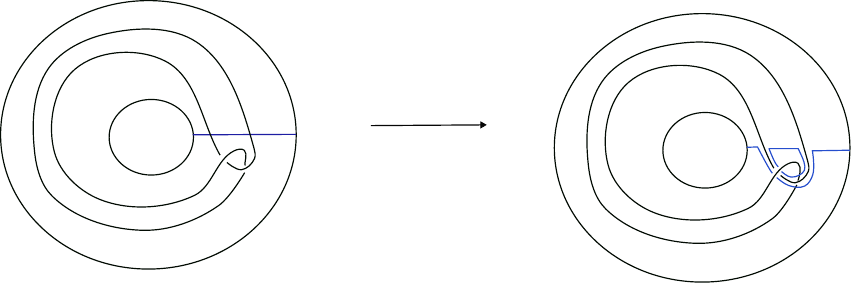}
\end{overpic}
\caption{Schematic of tubing on a whitehead clasp knot in a trivial bundle. The section to which we add a one handle is in blue.}\label{tubingpic}
\end{center}
\end{figure}

We now analyze what happens to $R''$ under the $\psi_2$ embedding. Since connected incompressible and boundary incompressible surfaces in $N$ are vertical tori, vertical annuli, and sections of $S$, the Thurston norm of a class in $H_2(N, \partial N; \bR)$ is realized by unions of such surfaces. By Poincar\'e-Lefschetz duality, there is a perfect pairing between $H_1(N; \bR)$ and  $H_2(N, \partial N; \bR)$ via algebraic intersection number. Recall that $t \in \pi_1(N)$ represents the fibre of $N$. We realize $t$ on some component of $\partial N$. Since $\left.(\psi_2^{-1} \circ \psi_1)\right|_{T_1 \cup \cdots \cup T_k}$ is (oriented) fibre preserving, it follows that $\psi_2^{-1}(t) = \psi_1^{-1}(t)$ and $[\psi_i(R'')] \cap [t] = \pm 1$ in the homology of $N$. Since only section-surfaces have non-trivial intersection with $[t]$, we see that $[\psi_i(R'')]$ is homologous to one section-surface and some number of vertical tori and annuli, i.e. the Thurston norm of $[\psi_i(R'')]$ is $|\chi(S)|$ for \emph{both} $i = 1, 2$. Thus, there are two subcases to consider $|\chi(R'')| > |\chi(S)|$ and $|\chi(R'')| = |\chi(S)|$.

{\bf Subcase $|\chi(R'')| > |\chi(S)|$.} This means that $x([R'']) > x([\psi_i(R'')])$ for both $i = 1, 2$. By \cite[Corollary 2.4]{Gabai}, there is at most one slope where the Thurston norm of an embedded surface can go down under filling. Thus, it follows that $s_1 = s_2$.

{\bf Subcase $|\chi(R'')| = |\chi(S)|$.} In this case, $R_i = \psi_i(R'')$ must be a horizontal section by the analysis above. Further, as we remarked, it follows that $R_i$ does not intersect $\eta_i$ for both $i = 1, 2$. Since $N = S \times \bS^1$, there is a fibre-preserving and orientation-preserving diffeomorphism $\varphi \in \Diffeo^+(N)$ that takes $R_1$ to $R_2$. Replacing $\psi_1$ with $\varphi \circ \psi_1$, we can therefore assume that a neighborhood of $R_1$ is fixed point-wise under $\psi_2 \circ \psi_1^{-1}$. Note, we are using the fact that $g$, and therefore  $\psi_2 \circ \psi_1^{-1}$, is orientation preserving here. After cutting $M$ along $R''$, we get two diffeomorphic knots complements in $S \times I$ where the diffeomorphism $\tau$ is the identity on $\partial (S\times I)$. This is because  $\tau$ is the identity on $S\times\set{0,1}$ and it is fibre-preserving on $\partial S\times I$ (and so isotopic to the identity). Moreover, neither $\eta_1$ or $\eta_2$ is contained in a $3$-ball in $S\times I$ by assumption. 

Since $\eta_1$ is not an unknot in $S \times I$, we can pick an embedding $\rho: S \times I \to \bS^3$ such that the image $\rho(\eta_1)$ is not an unknot in $\bS^3$. To choose such an embedding, observe that $S \times I$ is a handle body, so we can pick a maximal system of compressing disks and dual 1-handles. Since $\eta_1$ is not contained in a ball, it must intersect some compressing disk nontrivially with corresponding handle $h$. We may build $\rho$ by starting with a standard, unknotted embedding for all handles except for $h$, and then knot $h$ in a non-trivial way. By doing this, we can guarantee that $\rho(\eta)$ is a non-trivial satellite knot in $\bS^3$.
Using $\rho$, we can then extend $\tau$ to a diffeomorphism of $\bS^3$ taking $\rho(\eta_1)$ to $\rho(\eta_2)$.
By the Gordon-Luecke Theorem, it follows that $s_1 = s_2$.

From now on, we can assume that $m_1 \neq 0$ and $m_2 \neq 0$. 
 
\textbf{(2) Case $[\al_1]_{rel} \neq 0$ or $[\al_2]_{rel} \neq 0$ in $H_1(S, \partial S; \bZ)$}. 

In this case, we adapt the argument of Rodr\'igues-Migueles \cite{RM2020}. We may assume $[\al_1]_{rel} \neq 0$. The homology condition implies that there is an essential simple closed curve $\gamma \in S$ such that the algebraic intersection number $i(\al_1, \gamma) = n > 0$, see for example \cite{FM2011}. 

In $M$, this gives a punctured torus $Q = \psi_1^{-1}(\gamma \times \bS^1)$ whose boundary is homologous to $n$-times the meridian $\mu_1$ of $\eta_1$, where $n > 0$ by the previous paragraph. In particular, $[\mu_1]$ is torsion in homology. By half-lives and half-dies with $\bR$ coefficients, the homology kernel of the inclusion $\partial M \to M$ is Legendrian, so at most half of the homology of $T$ can die in $M$. That is, since any two linearly independent classes in $H_1(T; \bR)$ intersect non-trivially, the Legendrian property implies that $\ker(i_* : H_1(T; \bR) \to H_1(M; \bR))$ is either zero or one dimensional. Moreover, since $[\mu_1]$ is torsion under the inclusion $\iota: T\hookrightarrow M$, it must be in the kernel and therefore generates $\ker(i_* : H_1(T; \bR) \to H_1(M; \bR))$. If we also have that $[\al_2]_{rel} \neq 0$, then the meridians both generate this kernel and so $\mu_1 = \pm \mu_2$ and we are done. To finish this case, assume $[\al_1]_{rel} \neq 0$ and $[\al_2]_{rel} = 0$.
Then $\mu_1 =  p \mu_2 + q \ell_2$ for some $q \neq 0$ and we have that $[\partial \psi_2(Q)]$ in $H_1(N; \bZ)$ is $n q [\eta_2]$. Since $H_1(N; \bZ)$ is torsion free, we get that $[\eta_2] = 0$ in $H_1(N; \bZ)$. But this implies $m_2 = 0$, which is covered by the Case (1) and so again it follows that $s_1 = s_2$.

\textbf{(3) Case $[\al_1]_{rel} = 0$ and $[\al_2]_{rel} = 0$ in $H_1(S, \partial S; \bZ)$}. 

As before, let $T$ be the boundary torus that we are filling and let $T_1 \cup \cdots \cup T_k$ be the other boundary components of $M$. Let $R$ be a section of $N$ and pick a basis $\langle t, \rho_j \rangle$ for $\pi_1(T_j)$ with $\rho_j$ given by $\partial R$ and $t$ isotopic to a fibre. For each knot, we also have our basis $\langle \ell_i, \mu_i\rangle$ for $\pi_1(T)$ as above. By assumption, we must have that $\eta_i\cap R \neq \emptyset$ for each $i$.

Up to isotopy, we have that each component of $R \cap \partial \mathcal{N}(\eta_i)$ is a copy of $\pm \mu_i$, where $\mathcal{N}(\eta_i)$ is a tubular neighbourhood of $\eta_i$. Recall that $\eta_i = \al_i t^{m_i}$ and $[\al_i]_{rel} = 0$ in $H_1(S, \partial S; \Z)$. In particular, replacing $S$ by $R$, it follows that $[\al_j] = \sum_{j = 1} ^k a_{i,j} [\rho_j]$ in $H_1(R; \Z)$ for some $a_{i,j} \in \bR$ for $i = 1,2$. By homology considerations, we get the following equalities in $H_1(M; \bR)$:
\begin{enumerate}
\item[{\it (i)}] $m_i[\mu_i] = \sum_{j = 1}^k [\psi_i^{-1} \rho_j]$ -- by using $R$
\item[{\it (ii)}] $[\ell_i] = m_i [\psi_i^{-1}t] + \delta_i [\mu_i] + \sum_{j = 1} ^k a_{i,j} [\psi_i^{-1}\rho_j]$ -- by the above observation.
\end{enumerate}

Notice that $[t]$ is not a linear combination of $\{[\rho_j]\}_{j =1}^k$ in $H_1(N; \Z)$, so this is also true for their pullbacks in $H_1(M ;\Z)$ by both $\psi_1$ and $\psi_2$. Thus, since $m_i \neq 0$, we get that $[\mu_i]$ and $[\ell_i]$ are linearly independent in $H_1(M; \bR)$. Further, by hypothesis, $\left.(\psi_2^{-1} \circ \psi_1)\right|_{T_1 \cup \cdots \cup T_k}$ is the identity, so we have that $[\psi_1^{-1} \rho_j] = [\psi_2^{-1} \rho_j]$, which implies $m_2 [\mu_2] = m_1 [\mu_1]$. Let $\mu_1 = p \mu_2 + q \ell_2$ in $H_1(T; \bZ)$, then the first equation (i) becomes:

\begin{enumerate}
\item[{\it ($i^*$)}] $m_2 [\mu_2] = m_1 [\mu_1] = m_1(p [\mu_2] + q [\ell_2])$.
\end{enumerate}

Since $[\mu_2]$ and $[\ell_2]$ are linearly independent in $H_1(M;\Z)$, we get that $q = 0$. This gives $p = \pm 1$ since $\mu_i$ are both simple. Thus, $[\mu_1]= \pm [ \mu_2] \in H_1(T; \bZ)$, so the slopes are equal. 

Since we covered all possible cases, this concludes the proof of Theorem \ref{SF}.  \end{proof}

\brem[Condition (2) in Theorem \ref{SF}.]  The isotopy condition is necessary given the following example. Let $M_0=S_{g,2}\times \mathbb S^1$ and consider the slopes $(1,0)$ and $(1,q)$, $q\neq 0$, on one of the two boundary component $T_0$. Then, $M_0(1,0)$ and $M_0(1,q)$ are both trivial circle bundles and so there is a homeomorphism $f:M_0(1,0)\rar M_0(1,q)$. However, $f$ on the boundary is mapping the meridian $\mu$, the boundary of a section, to the section of $M_0(1,q)$, which twists $\mu$ $q$-times along the fiber direction of $M_0$. 

Take $M=M_0(1,0) \cong^+ M_0(1,q)$ and $N = M_0(1,0)(1,0) \cong^+ M_0(1,q)(1,-q)$. This gives two embeddings of $M$ in $N$ that differ by pre-composing by $f$ and the corresponding slopes are clearly different. This is precisely because $f$ does not satisfy condition $(2)$. It seems likely that if the the knots are not torus-satellites of a fibre or certain 1-bridge knots, then Theorem \ref{SF} would hold without condition (2).
\erem

\brem[Alternate approach to in Theorem \ref{SF}.] An alternative approach to prove Theorem \ref{SF} may be to use \cite[Corollary 2.14]{Gabai}, which states that a knot (not contained in a 3-ball) in a connected sum of $\mathbb{S}^2 \times \mathbb{S}^1$'s is determined by its complement. In the setting of Theorem \ref{SF}, one could fill the boundaries of $N$ along the fibers to obtain such a connected sum. However, even though the 3-ball restriction is resolved by Gordon and Luecke \cite{GL1989}, the issue of knots becoming trivial after filling or isotopic only after filling, requires delicate case analysis. We expect that a proof along these lines would have similar complexity to the proof given above.
\erem

To prove Theorem \ref{SFone}, the basic idea will be to exploit the fact that circle bundles over hyperbolizable surfaces only fibre in one way. We will show that the pre-images of some fibre from $M(s_i)$ in $M$ must coincide using group theoretic facts and then conclude that the slopes were equal using Theorem \ref{SF}. The centrality of the fibre in $\pi_1(M(s_i))$ will play a key role, so we begin with a Lemma about centers in one relator surface groups.

\blem\label{nocenter}  Let $S$ be an orientable surface with $\chi(S) < -1$ and $\al \in \pi_1(S)$. The group $G = \pi_1(S)/\llangle \al \rrangle$ has trivial center.
\elem

\bpf 

First, if $\al = 1$, then $G = \pi_1(S_1)$, which has trivial center as $\chi(S) < 0$. When $\al \neq 1$, we split this proof into the closed and non-closed cases.

{\bf Case 1}. $S$ is \emph{not closed}. 

Since $\pi_1(S)$ is a free group, $G$ is a one-relator group. By Murasugi \cite{M1964}, if $\pi_1(S)$ has at least three generators, then $G$ has trivial center. Since $\chi(S) < -1$, we know that it has at least 3 generators and we are done. Note, with only two generators Murasugi's result does not hold, which is why we need $\chi(S)<-1$.

{\bf Case 2}. $S$ is \emph{closed}. 

Here, we use a trick of Hempel \cite{H1990} as reformulated by Howie \cite{H2004} to reduce this to the previous case. We begin with an outline of Hempel's construction. 

First, Hempel shows that there is a non-separating simple closed curve $\beta$ such that the homological intersection number with $\alpha$ is trivial. This allows one to build an infinite cyclic cover $\widetilde{S} \to S$ where $\alpha$ lifts to a curve $\alpha'$. The cover $\widetilde{S}$ is obtained by cutting along $\beta$ and gluing copies of $S' = S \smallsetminus \beta$. Indexing these copies by $\mathbb{Z}$, we get a surface $F = S_0' \cup \cdots S_n'$ which minimally contains $\alpha'$, i.e. the lift $\alpha'$ is not contained in the surfaces $F_0 = S_0' \cup \cdots S_{n-1}'$ or  $F_1 = S_1' \cup \cdots S_{n}'$. Hempel then shows that $\pi_1(F_0)$ and $\pi_1(F_1)$ inject into $\pi_1(F)/\llangle \alpha' \rrangle$ and that $G$ is a HNN extension of $\pi_1(F)/\llangle \alpha' \rrangle$ with associated subgroups $\pi_1(F_0)$ and $\pi_1(F_1)$ (more technically they are embeddings, but we will abuse notation here). 
Fixing a presentation $\pi_1(F)/\llangle \alpha' \rrangle = \langle A \mid R \rangle$, the HNN extension gives the presentation 
\[G = \left\langle A, t \mid R, t h t^{-1} = \phi(h) \; \forall h \in \pi_1(F_0)\right\rangle,\]
where $\phi$ is the isomorphism $\phi: \pi_1(F_0) \to \pi_1(F_1)$ obtained from translating in the cyclic cover.

Since $\pi_1(F)$ is free and $\chi(F) < -1$, we know that $\pi_1(F)/\llangle \alpha' \rrangle$ has trivial center by Case 1. It remains to show that this HNN extension also has trivial center. For this, we use Britton's Lemma \cite{LS1977}. 

First, notice that any element of an HNN extension can be written as 
\[w = g_0 t^{s_0} \cdots g_n t^{s_n}g_{n+1}\]
 with $s_i \in \{\pm 1\}$ and  $g_i \in \pi_1(F)/\llangle \alpha' \rrangle$. The term $g_i = e$ is allowed but there are no terms of the form $tht^{-1}$ for $h \in  \pi_1(F_0)$ or $t^{-1}ht$ for $h \in  \pi_1(F_1)$. Consider a non-trivial $w \in Z(G)$ and assume that in the presentation above $n$ is minimal. By Britton's Lemma, as $w$ is non-trivial, it follows that $n > 0$. Since $w$ is in the centre, we can assume that $g_{n+1} = e$ as we may cyclically move $g_{n+1}$ to the front and absorb it into $g_0$. Our next trick is to push any $g_i \in \pi_1(F_0) \cup \pi_1(F_1)$ for $i > 0$ ``past'' $t$ terms until we cannot anymore. Abusing notation for the new $g's$, this gives $w$ as a product of terms $g_it^{s_i}$ where $g_i = e$ or $g_i \notin \pi_1(F_0) \cup \pi_1(F_1)$ for $i > 0$. Absorbing any terms that are pushed out on the left or right into $g_0$, we are left with two cases:

\begin{enumerate}
\item $w = g_0 t^r$ for some $g_0 \in \pi_1(F_0) \cup \pi_1(F_1)$ and some $r$
\item  $w$, after maybe cyclically permuting, has the same form as before but $g_0 \notin \pi_1(F_0) \cup \pi_1(F_1)$ and $g_i \notin \pi_1(F_0) \cup \pi_1(F_1)$ or $g_i = e$ for $i > 0$.  
\end{enumerate}

In the first case (1), since $w$ is central, we would have $e = [w,t] = g_0 t^r t t^{-r} g_0^{-1} t^{-1} = g_0tg_0^{-1}t^{-1}$. However, since $t$ acts as a shift map $\pi_1(F_0) \to \pi_1(F_1)$, we know that it does not commute with any non-trivial $g_0$. If $g_0$ is trivial, then $w = t^q$, but for this to be in the center we must have $q = 0$, contradicting non-triviality of $w$.

Now we focus on the latter case (2). Consider 
\[e = [w,t] = g_0 t^{s_0} \cdots t^{s_{n-1}} g_n t g_n^{-1}t^{-s_{n-1}}\cdots g_1^{-1} t^{-s_0}g_0^{-1} t^{-1},\]
where $g_i \notin \pi_1(F_0) \cup \pi_1(F_1)$ or $g_i = e$ for $i > 0$. Britton's Lemma implies that one of the factors must have a ``bad term'' of the from $tht^{-1}$ for $h \in  \pi_1(F_0)$ or $t^{-1}ht$ for $h \in  \pi_1(F_1)$. In the sub-case where $g_n \notin \pi_1(F_0) \cup \pi_1(F_1)$, minimality of $n$ and the fact that none of the original $s_i$ and $g_i$ give a ``bad term'', we get a contradiction. Note, the pushing operation cannot introduce ``bad terms,'' since we assume minimality of $n$. In the case where $g_n = e$, we obtain a cancellation and must have that $g_{n-1}$ is part of a ``bad term.'' But this is again only possible if $g_{n-1} = e$, leading to further cancellation. Repeating this process gives $[w,t] = g_0tg_0^{-1}t^{-1}$ and we can conclude as before.

Obtaining a contradiction in both cases, we conclude that the center of $G$ is trivial. \epf

We can reinterpret Lemma \ref{nocenter} in the context of circle bundles as follows.

\bcor\label{commutator}   Let $S$ be an orientable surface with $\chi(S) < -1$ and let $N$ be an orientable circle bundle over $S$. Fix $\beta \in \pi_1(N)$. If $\gamma \in \pi_1(N)$ has the property that the commutator $[g, \gamma] \in \llangle \beta \rrangle$ for all $g \in \pi_1(N)$, then $\gamma  \in  t^r  \llangle \beta \rrangle$ where $t$ generates the center of $\pi_1(N)$.
\ecor

\bpf Let $C= Z(\pi_1(N)) = \langle t \rangle$ and $H = \llangle \beta, t \rrangle$. Consider the quotients $\pi_1(S) = \pi_1(N)/C$ and $G = \pi_1(N)/H =  \pi_1(S)/ \llangle \beta C \rrangle$. From our assumption, it follows that $[gH, \gamma H] = 1$ for all $gH \in G$. In particular, $\gamma H \in Z(G)$. By Lemma \ref{nocenter}, it follows that $\gamma H = 1$, so $\gamma \in H = \llangle \beta, t \rrangle$. Since $t$ is central in $\pi_1(N)$, $\gamma \in t^r \llangle \beta \rrangle$ for some $r$.
\epf

Our remaining necessary group theoretic Lemma has to do with words being trivial in normal closures of elements in surface groups.

\blem\label{powersum} Let $S$ be an orientable surface and $w \in \pi_1(S)$. Write $\beta \in \llangle w \rrangle$ as \begin{equation}\label{comeq} \beta = \prod_j g_j^{-1} w^{\eps_j} g_j \quad \text{ for } \quad g_j \in \pi_1(S) \text{ and } \eps_j = \pm 1.\end{equation}
If $\beta = 1$ and $w \neq 1$, then $\sum_j \eps_j = 0$.
\elem

\bpf If $[w] \in H_1(S; \mathbb{Z})$ is non-zero, then we are done as $[\beta] = \left(\sum_j \eps_j\right) [w]$. Thus, we need to deal with the case where $w$ is null-homologous. 

To do this, consider the following tower of commutators, i.e. derived series of $\pi_1(S)$. Let $H_1 = [\pi_1(S), \pi_1(S)]$ and $H_{i+1} = [H_i, H_i]$ for $i > 1$. Since the commutator subgroup $H_1$ of $\pi_1(S)$ is free (of infinite-type), we know that this tower is never trivial. Moreover, $\cap_{j=1}^\infty H_j = \{e\}$, i.e. the perfect core of a surface group is trivial. Thus, since $w \in H_1$ and $w \neq 1$, there is a $p$ at which $w \in H_{p}$ but $w \notin H_{p+1}$. Note, since $H_p$ is normal, we know that $g_j^{-1} w^{\eps_i} g_j \in H_p$ as well. Let $S_p$ be the cover of $S$ corresponding to $H_p$. Since $w \in H_1$, we will make the following convenient change: if $S$ is closed, we will replace $S$ by $S_1$ and $w$ by a lift. This allows us to assume that $\pi_1(S)$ is free. Note, under this change, Equation \ref{comeq} is no longer a product of just commutators, but also includes translates of $w$ by $g_j$ that act like deck transformations.

Consider the deck group $G_p = \pi_1(S)/H_p$. We know that $w$ and its conjugates lift to $S_p$ and that the homology class $[w]_p \in H_p/H_{p+1} = H_1(S_p; \bZ)$ is nontrivial\footnote{We are using $[x]_p$ to underline the fact that we are considering the homology in $H_1(S_p;\Z)$.}. To see that $\sum_j \eps_j = 0$, we switch from $S$ to a CW-complex. Since $\pi_1(S)$ is free, we can assume $S$ is just a wedge of $n$ orientated loop $\{e_l\}_{l = 1}^n$ with a unique basepoint $\star$. We allow $n = \infty$ and we assume that the loops form a free basis for $\pi_1(S)$. 

Lifting to $S_p$, we fix a lift $\wt{\star}$ and see that the $0$-chains $C_0(S_p; \bZ) = \mathbb{Z}(G_p) \cdot \wt{\star}$ is a 1-dimensional $\mathbb{Z}(G_p)$-module. Similarly, since $G_p$ acts freely and properly discontinuously, the orbits of the orientated lifts $\wt{e}_l$ of $e_l$ based at $\wt{\star}$ are disjoint, so $C_1(S_p; \mathbb{Z}) = \mathbb{Z}(G_p)^n$. Realizing the lift $\wt{w}$ of $w$ based at $\wt{\star}$ as a 1-chain, we can write \[\wt{w} = \sum_{l = 1}^n \left(\sum_{i =1}^{m_l}  \alpha_{l,i}  h_{l,i} \right) \cdot \wt{e}_l,\] were $h_{l,i} \in G_p$ are distinct for fixed $l$, $m_l \geq 1$ for finitely many $l$, and $a_{l,i} \in \bZ$ are nonzero. By assumption, $[\wt{w}]_p$ is a cycle and it is nonzero in $H_1(S_p; \bZ)$. Abusing notation, we will let $g_j$ denote the image of $g_j$ in $G_p = \pi_1(S)/H_p$. Then, Equation \ref{comeq} with $\beta = 1$ and the fact that $C_2(S_p; \bZ) = 0$, implies that the chain  \[ \sum_j \eps_j g_j \sum_{l=1}^n \left(\sum_{i = 1}^{m_l} \alpha_{l,i}  h_{l,i} \right) \cdot \wt{e}_l  = 0\]

We can now group the $g_i$'s into a collection $\{g_k'\}_{k = 1}^q$ of distinct elements and group the $\eps_j$ coefficients into $\beta_k = \sum_j \eps_j \delta(g_j, g_k')$ where $\delta(u,v) = 1$ if $u = v$ and $0$ otherwise. This gives \[ \sum_{l=1}^n \left(\sum_{k=1}^q \sum_{i=1}^{m_l} \alpha_{l,i} \beta_k g_k' h_{l,i} \right) \cdot \wt{e}_l  = 0.\]
If $\beta_k = 0$ for some $k$, then we can drop $g_k'$ form the list. Thus, we may assume $\beta_k \neq 0$ for all $k$.

Notice that $G_p$ is a (tower of) semidirect product(s) of torsion-free abelian groups. Since each such group is left-orderable, so is $G_p$. In particular, we may assume that $h_{l,1} < h_{l, 2} < \cdots < h_{l, m_l}$ for each $l$. For fixed $l$ and $k$, it follows that amongst $g_k' h_{l,i}$, the maximal element is $g_k' h_{l, m_l}$. For fixed $l$, let $s_l$ be such that $g_{s_l}'h_{l, m_l}$ is maximal amongst all $g_k' h_{l,i}$. Note, $s_l$ is unique since $g_k' h_{l, m_l} =  g'_t h_{l, m_l}$ implies $g_k' = g_t'$ and $k  = t$.

In the sum above, each coefficient of $\wt{e}_l$ must be zero. In particular, by uniqueness of $s_l$, we have that $\alpha_{l,m_l} \beta_{s_l} g_{s_l}' h_{l,m_l}  = 0$. However, this is a contradiction since $\alpha_{l,m_l} \beta_q \neq 0$ by assumption. Thus, all $\beta_k = 0$ and so $\sum_j \eps_j = \sum_k \beta_k = 0$. 
\epf

We now turn to showing that Theorem \ref{SFone} reduces to Theorem \ref{SF}. Our main step will be to show Proposition \ref{fibre}, which states that up to the choice of pre-image, the vertical fibers of the filled manifold pulled back under $\psi_1$ and $\psi_2$ are isotopic in $M$. After drilling this pre-image of a fibre, we will be able to apply Theorem \ref{SF}.\bprop\label{fibre} Let $M$ be a compact orientable 3-manifold where $\partial M$ is a union of tori.
Let $s_1$ and $s_2$ be slopes on the first boundary component of $M$. Assume that
\begin{enumerate}
 \item $M(s_i)$ is an orientable circle bundle over an orientable surface $S_i$ for $i = 1, 2$,
 \item $\chi(S_1) < -1$.
 \end{enumerate}
 Then there is a loop $\gamma$ in $M$ whose image in $M(s_1)$ and $M(s_2)$ is homotopic to a fibre.

\eprop

\begin{proof}
 Let $N_i = M(s_i)$ for $i = 1,2$ and consider the embeddings $\psi_i:M \to N_i$ realizing $M$ as the complement of an open tubular neighborhood of the cores $\eta_i$ of the filling, for $i = 1,2$ respectively. We can assume $s_1 \neq s_2$ as otherwise the conclusion is trivial.

Let $T$ be the torus boundary that we are filling. Place a basepoint on $T$ and realize $s_i$ as $\sigma_i \in \pi_1(M)$ on $T$. For $i = 1, 2$, fix a circle bundle structure on $N_i$ and let $\gamma_i$ be a representative of this fiber in $\pi_1(M)$, pulling back via the embedding $\psi_i$. Notice that the choice of $\gamma_i$ is only well defined up to taking products of $\gamma_i$ with conjugates of $\sigma_i$.

Our goal will be to show that, up to the choice of pullback, $\gamma_1 = \gamma_2$ in $\pi_1(M)$.
One of the key tools of this argument will be Corollary \ref{commutator}. Further, by Lemma \ref{type}, we have that $S_1$ and $S_2$ have the same topological type. In particular, $\chi(S_1) = \chi(S_2) < -1$.

Let $H_i = \llangle \sigma_i \rrangle$ be the normal closure of $\sigma_i$ in $\pi_1(M)$. Observe that Dehn filling along $s_i$ corresponds to $\pi_1(N_i) = \pi_1(M)/H_i$.
Since both $N_1$ and $N_2$ are circle bundles, we know that $\gamma_i H_i \in \pi_1(N_i)$ is the element that uniquely (up to sign) generates the centre of $\pi_1(N_i)$, for each $i= 1,2$ respectively. This implies the following two key facts:
\begin{enumerate}
\item first, our choice of $\gamma_i$ is equivalent to picking elements from the coset  $\gamma_i H_i$;
\item\label{com} second, the centrality of $\gamma_i H_i$ in $\pi_1(N_i)$, means that in $\pi_1(M)$ the commutator $[g, \gamma_i] \in H_i$ for all $g \in \pi_1(M)$.
\end{enumerate}

Looking now at $\pi_1(N_1)$, fact (\ref{com}) means that $[gH_1, \gamma_2 H_1] \in \llangle \sigma_2 H_1 \rrangle$ for all $g \in \pi_1(M)$. 
By Corollary \ref{commutator}, we must have that $ \gamma_2 H_1 \in  \gamma_1^{r_1} \llangle \sigma_2 H_1 \rrangle$ for some $r_1$. Applying the same argument in $\pi_1(N_2)$, we get that $\gamma_1H_2 \in  \gamma_2^{r_2} \llangle \sigma_1 H_2 \rrangle$ for some $r_2$. In particular, there are $u_1, u_2 \in H_1$ and $v_1, v_2 \in H_2$ such that:
\begin{equation}\label{gammas}\tag{5.1}\gamma_2 = \gamma_1^{r_1} v_1 u_1 \quad \text{ and } \quad \gamma_1 = \gamma_2^{r_2} u_2 v_2\quad\text{in }\pi_1(M).\end{equation}
Note, we are using normality of $H_1$ and $H_2$ to move elements across each other and write these expressions. It follows that $\gamma_1 = \gamma_1^{r_1 r_2} v_3u_3$ and $\gamma_2 = \gamma_2^{r_1r_2} u_4 v_4$ for some $u_3, u_4 \in H_1$ and $v_3, v_4 \in H_2$, again using normality $H_1$ and $H_2$.

{\bf Claim 1.} $r_1 = r_2 = \pm 1$ or $\eta_1$ is homotopic to a non-zero power of the fibre.

\bpfc Looking in $N_1$, we have that$\gamma_1 H_1 = \gamma_1^{r_1 r_2} \beta$ for some $\beta \in \llangle \sigma_2 H_1 \rrangle$. Let $f_1: \pi_1(N_1) \to S_1$ be the fiberwise projection. Applying it to both sides, we see that $1 = f_1^*(\beta)$, since the push-forward $f_1^*$ kills $\gamma_1H_1$. Recall that there is a presentation: 

\small
\[\pi_1(N_1) = \langle a_1, b_1, \ldots, a_{g}, b_{g}, d_1, \ldots d_k, t \mid a_i t = t a_i,  \, b_i t = tb_i, \,  d_it = td_i, \, \Pi[a_i, b_i]  \Pi d_i = t^{e_1} \rangle,\]
\normalsize

where $e_1$ is the Euler number of $N_1$. Note, that in our context, $t = \gamma_1 H_1$ is the fibre in $\pi_1(N_1)$.  This means that we can write $\sigma_2 H_1 = wt^q$ for some integer $q$ and {\it surface word} $w$ (i.e. $w$ is a word in only $a_i$'s, $b_i's$ and $d_i$'s). In this context, since $\beta \in \llangle \sigma_2 H_1 \rrangle=\llangle wt^q \rrangle$, we have:

\[\beta = \prod_j g_j w^{\eps_j} g_j^{-1} t^{q \sum_j \eps_j} \quad \text{ for some surface words } g_j \text{ and } \eps_j = \pm 1.\]

Notice that $1=f_1^*(\beta) = \prod_j g_j w^{\eps_j} g_j^{-1}$, where $f_1^*$ of a surface word is just itself. If $w \neq 1$, then by Lemma \ref{powersum}, we must have that $\sum_j \eps_j = 0$, which implies that $\beta$ is a surface word. Since it has trivial projection, then $\beta$ is also trivial in $\pi_1(N_1)$. Applying this to $\gamma_1H_1$, we get that $\gamma_1H_1 = \gamma_1^{r_1 r_2}H_1$, so $r_1 = r_2 = \pm 1$ as $\pi_1(N_1)$ has no torsion.
 
Assume now that $w = 1$, so $\sigma_2 H_1 = t^q$. If $q = 0$, we once again get $\beta$ is trivial in $\pi_1(N_1)$, so $r_1 = r_2 = \pm 1$.
 
Notice that under the $s_1$ Dehn filling, the curve $\sig_2$ maps to some power of the core $\eta_1$. That is, $\eta_1^k = \sig_2H_1$ in $\pi_1(N_1)$ for some $k$ with $k \neq 0$ since we have assumed $s_1 \neq s_2$. Thus, in our remaining case where $q \neq 0$ and $\sig_2H_1 = t^q$, we see that $\eta_1$ must be homotopic to some non-zero power of the fibre. By symmetry, we may assume that $\eta_2$ is a non-zero power of the fibre in $N_2$.
\epfc

We now deal with the case in which $\eta_1\simeq t^s$ in $\pi_1(N_1)$ and similarly, by symmetry, $\eta_2\simeq t^v$ in $\pi_1(N_2)$.
Note that the shadow $\hat \eta_i$, $i=1,2$, are both null-homotopic in the base surface of the bundle. 

Let $F$ be the essential sub-surface of $S_1$, the base surface of $N_1$, filled by the shadow of the $\eta_1$ and let $\Sigma = F^C$.

\textbf{Case 1.} Assume that some component of $\Sigma$ has negative Euler characteristic.

Let $W$ be the trivial bundle over such a component of $\Sigma$ with boundary $V$, which is separating in $N_1$. Since the fibration on $W$ is unique, $W$ has the following property: the only non-trivial simple closed curve on one component of $V$ that is isotopic into another component is a fibre. Further, $V$ is a separating union of vertical tori, so $\psi_2 \circ \psi_1^{-1}(V)$ is $N_2$ is also a separating union of vertical tori. Realizing $M$ as $W \cup_V (N_1 \smallsetminus  (\mathcal{N}(\gamma_1) \cup W))$, we see that after filling along $s_2$, the induced fibration of the pair $V \subset W$ has the desired isotopy property. In particular, the fibrations induces by $\psi_1$ and $\psi_2$ must be the same on $W$ and so a there is a simple closed curve in $M$ mapping for a fibre under each filling.

\textbf{Case 2.} Assume every component of $\Sigma$ has non-negative Euler characteristic.

In this case, $F$ is obtained from $S$ by removing a (possibly empty) collection of disks and annuli from $S$ and so it has negative Euler characteristic.  Consider a vertical torus $V_\gamma$ above a simple closed curve $\gamma\subset S_1$. By assumption on $\eta_1$, we have that the algebraic intersection number of $V_\gamma$ and $\eta_1$ is zero. Let $V_\gamma'$ be an embedded closed surface obtained from $\psi_1^{-1}(V_\gamma)$ in $M$ by attaching tubes (i.e. annuli) to $\partial \psi_1^{-1}(V_\gamma)$ in a regular neighborhood of an annulus containing $\partial \psi_1^{-1}(V_\gamma) \subset \partial M$. After the $s_1$ filling, we know that $[\psi_1(V_\gamma')] = [V_\gamma]$ in $H_2(N_1)$. We claim that by our assumption, there is a $\gamma$ such that $x([V_\gamma']) > 0$ in $H_2(M, \partial M \setminus T)$ and so the Thurston norm goes down under the $s_1$ filling. Indeed, assume $x([V_\gamma']) = 0$ and let $R$ be a norm-minimizing surface. Since $M$ is irreducible by our assumption on $\eta_1$, we can assume $R$ is a union of tori, which after the filling either compress or are vertical in $N_1$. Since every component of $\Sigma$ is either a disk or an annulus and these tori would be disjoint from $T$, after filling, the vertical pieces must be above annular components of $\Sigma$ and the compressible ones are either parallel to $\eta_1$ or above disk components of $\Sigma$. We can pick a simple closed curve $\beta$ on $F$ and pick $\gamma$ on $S$ such that $i(\beta, \gamma) = 1$ and $[\gamma] \neq 0$. This is possible since $F$ has negative Euler characteristic and $\Sigma$ is a union of disks and annuli. Lifting $\beta$ to $N_1$, we see that $i([V_\gamma], [\beta]) = i([\psi_1(R)], [\beta])= 1$. This implies that $\psi_1(R)$ has no vertical torus components. Further, our choice of $\gamma$ also eliminates all compressible components since $[\psi_1(R)] = [V_\gamma]$ in $H_2(N_1, \partial N_1)$. This implies that $R$ is empty, which is not possible as $[R] \neq 0$. Thus, $x([V_\gamma']) > 0$ in $H_2(M, \partial M \setminus T)$.

Lastly, we argue that $\psi_2(R)$ must compress. Indeed, since $\eta_2$ is a  homotopic to a power of a fibre and is disjoint from $\psi_2(R)$, we know that $[\psi_2(R)] = 0$ in $H_2(N_2, \partial N_2)$.  By \cite[Corollary 2.4]{Gabai}, we conclude that $s_1 = s_2$ and so we can choose our fibre pre-image as desired.

With all the cases considered, we can assume that we have chosen homotopically equivalent $\gamma_1$ and $\gamma_2$ at the beginning as the pre-images of the fibre structure. \end{proof}

We can now prove Theorem \ref{SFone}.

\begin{proof}[Proof of Theorem \ref{SFone}.] Let $N = M(s_1) \cong^+ M(s_2)$ and let $\eta_i$ be the two knots. Recall that $H_1(N)=\Z^{2g}\oplus \quotient{\Z}{e\Z}$ for $e$ the Euler number of the bundle. Let $\mu_i, \sigma_i$ denote the slopes on $\partial M$ corresponding to the meridians and the chosen longitudes of $\eta_i$ as in Theorem \ref{SF}.

We will deal with various cases. We first rule our knots that are non-torsion in homology and then use tropological and algebraic augments for the other cases.

\textbf{Case 1.} The class $[\eta_1]$ or $[\eta_2]$ is \emph{not} torsion in $H_1(N)$. 

For simplicity, we will abuse notation and use loops and their corresponding classes in homology interchangeably. We are free to assume $\eta_1$ is not torsion in $H_1(N)$ and therefore $\sigma_1$ is not torsion in $H_1(M)$. By half-lives half-dies in $M$, it follows that $\mu_1$ is torsion in $H_1(M)$. That is, there exists $\omega_1 \in H_2(M,\partial M)$ with $\partial\omega_1=k\mu_1$. Mayer-Vietoris on $M\cup_T V=N$, for $V$ a solid torus with $T = \partial V$, gives:
\[ H_2(N) \overset{ \delta_2}\longrightarrow H_1(T)  \overset{ \iota}\longrightarrow H_1(V)\oplus H_1(M)  \overset{\xi}\longrightarrow H_1(N)  \overset{ \delta_1}\longrightarrow H_0(T),\]
which becomes:
\[H_2(N)  \overset{\xi}\longrightarrow 
 \langle \sigma_1\rangle  \oplus \langle \mu_1 \rangle  \overset{\iota}\longrightarrow
  \langle \sig_1 \rangle \oplus H_1(M)  \overset{\xi}\longrightarrow 
  \Z^{2g}\oplus \quotient{\Z}{e\Z} \longrightarrow  0,\]
where $\delta_1$ is the zero map since $T$ is separating. Moreover, since $\im(\delta_2)=\langle k\mu_1\rangle$ we can rewrite the long exact sequence as:
\[0 \longrightarrow  \langle \sigma_1\rangle \oplus \langle \mu_1\rangle / {k \langle \mu_1\rangle}  \overset{\iota}\longrightarrow  \langle \sig_1\rangle \oplus H_1(M)  \overset{\xi}\longrightarrow  \Z^{2g}\oplus \quotient{\Z}{e\Z} \longrightarrow 0,\]
and $\iota(0, \mu_1) = (0, \mu_1)$, $\iota(\sigma_1,0)=(\sigma_1, \sigma_1)$. Thus, we obtain that:
\[\quotient{H_1(M)}{\left(\quotient{\langle \mu_1\rangle}{k \langle \mu_1\rangle}\right)}\cong \Z^{2g}\oplus \quotient{\Z}{e\Z}\]
which gives one of the two options:
\[ H_1(M) \cong \Z^{2g}\oplus \quotient{\Z}{e\Z}\oplus \quotient{ \Z}{k\Z}\qquad H_1(M) \cong\Z^{2g}\oplus \quotient{ \Z}{ek \Z}.\]

In either case, any slope $s_1=p\mu_1+q\sigma_1$ with $p\neq 0$ and $\abs q >0$ is not torsion in $H_1(M)$. Hence, $\text{rank}(H_1(M(s_2)))<2g$ which implies that $M(s_2)$ cannot be homeomorphic to $N$ and so the slopes have to be both equal to $(1,0)$.

We now deal with the case where both knots are torsion (including zero) in $H_1(N)$. We will be interested in how $\eta_i$ interacts with non-separating tori $\Sigma$ in $N$.

\textbf{Case 2.} There is an embedded, incompressible, non-separating torus $\Sigma$ in $N$ that is disjoint from $\eta_1$. 

 Abusing notation, we consider $\Sigma$ as lying in $M$. Since every non-separating torus in $N$ is incompressible, $\psi_i(\Sigma)$ is incompressible in both $M(s_1)$ and $M(s_2)$ and so vertical. Let $\mathcal F_i$ be the induced foliation of $\psi_i(\Sigma)$ by the fibers of $N$. By Proposition \ref{fibre}, there is a simple loop $\gamma$ in $M$ whose image in $N$ is a homotopic to a fibre. In particular, $\psi_i(\gamma)$ is homologous in $N$ to a fibre on $\psi_i(\Sigma)$. Since both $\eta_1$ and $\eta_2$ are torsion in $H_1(N)$, it follows that $\gamma$ is is homologous to the fibers of $\mathcal F_i$ in $H_1(M; \mathbb{R})$ for both $i = 1$ and $i = 2$. This implies that the fibers of $\psi_1(\mathcal F_1)$ and $\psi_1(\mathcal F_2)$ are homologous in $H_1(M(s_1); \mathbb{R})$. However, since $N = M(s_1)$ is a circle bundle over a surface of negative Euler characteristic and $\psi_1(\Sigma)$ is a vertical torus, this is only possible if $\psi_1(\mathcal F_1)$ and $\psi_1(\mathcal F_2)$ coincide. We can therefore assume that $\gamma$ was simply a fibre of $\mathcal F_1$ and it follows that $\psi_i(\gamma)$ is isotopic to a fibre in $M(s_i)$.

We can assume that $\left.\psi_1\right|_\gamma  = \left.\psi_2\right|_\gamma$ after possibly applying a fibre-reversing (but orientation-preserving) diffeomorphism of $N$. Further, we can assume that there is a regular neighborhood $V$ of $\gamma$ such that $W = \psi_1(V) = \psi_2(V)$ is a fibered vertical torus in $N$. By construction, we must have that $\left.(\psi_2 \circ \psi_1^{-1})\right|_W$ sends fibers to fibers. Thus, in particular, $\left.(\psi_2^{-1} \circ \psi_1)\right|_{\partial W}$ must be the identity up to isotopy since the core fibre is fixed point-wise. This gives the final necessary condition to apply Theorem \ref{SF} to conclude that $s_1 = s_2$ or that the knots are both trivial.

It remains to deal with the case where both $\eta_1$ and $\eta_2$ and torsion in $H_1(N)$ and non-trivially intersect every embedded, incompressible, non-separating torus $\Sigma$ in $N$.

\textbf{Case 3.} Assume that for every embedded, incompressible, non-separating torus $\Sigma$ in $N$, the geometric intersection $\iota(\Sigma,\eta_1) > 0$. 

By Case 1, $\eta_1$ is torsion in $H_1(N)$ and so the shadow $\al_1$ of $\eta_1$ is zero $H_1(S)$. Since $\Sigma$ must be vertical in $N$, it follows that the algebraic intersection pairing of $\eta_1$ and $\Sigma$ is always zero.  We now fix a particular $\Sigma$.

Let $\bar \Sigma$ be the embedded closed surface in $M$ be obtained from $\Sigma \cap \psi_1(M)$ by attaching ``tubes'' (i.e. annuli) that fellow travel $\partial \mathcal N(\eta_1)$ along nested pairs of boundary components and then taking the pre-image in $M$. Similarly, let $\hat \Sigma$ be the embedded surface in $M$ obtained by ``tubing'' a Thurston-norm minimizing representative of $\psi_1^{-1}(\Sigma)$ in $H_1(M, \partial M; \mathbb{R})$, as was done in Case (1) of the proof of Theorem \ref{SF}. In particular, $\hat \Sigma$ is a Thurston-norm minimizing surface with $x(\hat\Sigma) = \abs{\chi (\hat \Sigma)}$ and it is homologous to $\bar \Sigma$ in $H_1(M, \partial M\; \mathbb{R})$. It follows that $\psi_1(\bar \Sigma)$ is homologous to $\psi_1(\hat \Sigma)$ in $H_1(N; \mathbb R)$. Further, since $\bar \Sigma$ is homologous to $\psi_1^{-1}(\Sigma)$ in $H_1(M, \partial M\; \mathbb{R})$, we see that  $\psi_1(\hat \Sigma)$ is homologous to $\Sigma$ in $H_1(N; \mathbb R)$.

Notice that our assumption guarantees that $x(\hat\Sigma) > 0$ in $M$, as otherwise $\hat\Sigma$ is a union of incompressible, non-separating tori in $M$, which would imply that $\psi_1(\hat\Sigma)$ is an incompressible, non-separating union of tori that misses $\eta_1$, a contradiction. It follows that the Thurston norm of $\hat\Sigma$ drops under the $s_1$ filling. We can see that the norm must also drop under the $s_2$ filling as follows. If $N_2$ has non-trivial Euler number, then the only impossible surfaces in $N_2$ are vertical tori, and so the Thurston norm of any class is zero. Alternatively, if $N_2$ is a trivial bundle, then the only possibility is that $\psi_2(\hat\Sigma)$ compresses to have section surfaces. But these are disjoint from $\eta_2$, and so the argument of Case (1) of Theorem \ref{SF} applies. Thus, by \cite[Corollary 2.4]{Gabai}, we conclude that $s_1 = s_2$.
\end{proof}

\section{Application to canonical links}\label{appsec}
In this section we want to prove the following application of Theorem \ref{SFone}:

\appl*

\subsection{Reidemeister moves for canonical diagrams}\label{reiddiffcan}

First, we focus on understanding how diagrams of canonical link behave under isotopy and diffeomorphism of the ambient tangent bundle.

{\bf Legendrian and loop Reidemeister moves}.

Let $N = UT(S)$ or $PT(S)$ as above and consider two canonical links $\hat\kappa_1$ and $\hat\kappa_2$ in $N$ corresponding to diagrams $\kappa_1$, $\kappa_2$ on $S$. Assume that there is $\psi \in \Mod(N)$ such that $\psi$ maps $\hat\kappa_1$ to $\hat\kappa_2$ up to isotopy. By Section \ref{bkg}, we can compose $\psi$ with (a lift of) a surface mapping class $\phi$ such that $\phi \circ \psi$ is a transvection. Thus, up to the action of $\Mod(S)$, we need to understand how (the diagrams) $\kappa_1$ and $\kappa_2$ are related given that  $\hat\kappa_1$ and $\hat\kappa_2$ are the same up to a transvection and isotopy.

\begin{figure}
\begin{center}
\begin{minipage}{0.4\textwidth}
\begin{overpic}[scale=1]{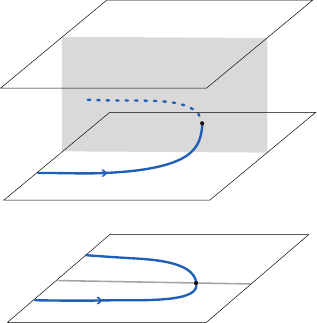}
\end{overpic}
\end{minipage}
\quad\quad
\begin{minipage}{0.4\textwidth}
\begin{overpic}[scale=1]{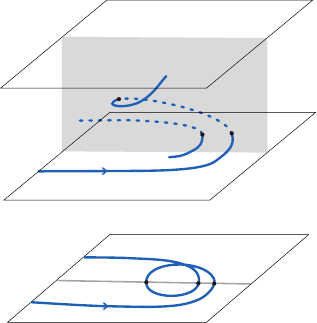}
\end{overpic}
\end{minipage}
\end{center}
\caption{A transvection adding a left loop. The right picture has an isotopy applied after a transvection about the grey vertical annulus/torus. The points are there to help visualize the intersections with the vertical annulus/torus.}
\label{tbtwist}
\end{figure}

Transvections act a little differently for $UT(S)$ and $PT(S)$. As before, pick a representative $\al$ of $[\al] \in H_1(S, \partial S; \bZ)$ as a weighted collection of disjoint simple closed curves and simple proper arcs $\al = \sum a_i \gamma_i$. The transvection corresponding to $[\al]$ corresponds to twisting $a_i$-times about the tori or annuli above each $\gamma_i$. When dealing with knots in $PT(S)$, at an intersection point of $\kappa_1$ with $\gamma_i$, every twist adds a left or right cusp. Thus, depending on how $\kappa_1$ meets the collection $\gamma_i$, we may be adding cusps by applying a transvection. In the case of $UT(S)$, instead of adding a cusp, a single twist adds a loop on the right or left, see Figure \ref{tbtwist}. This looks like a classical Reidemeister I move, but arises from an ambient diffeomorphism. Note, we will think of transvections as only adding loops or cusps. Removal will be handled by destabilizations introduced below.

\begin{figure}[htb!]
\centering
\begin{overpic}[scale=0.5]{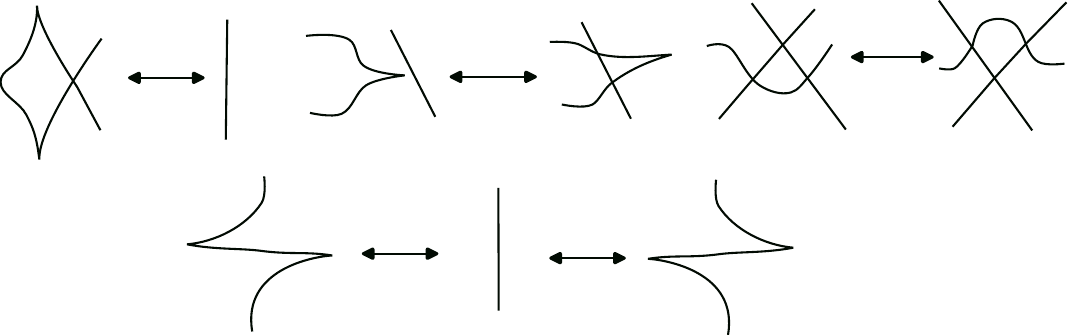}
\end{overpic}
\caption{ Legendrian Reidemaister moves and (de)stabilizations. }\label{lgmoves}
\end{figure}

It remains to understand the effect of an isotopy of $\hat\kappa_1$ in $N$ on the diagram $\kappa_1$. In the context of $PT(S)$, any two diagrams of canonical knots that are isotopic in $PT(S)$ must be equivalent up to Legendrian Reidemeister moves and Legendrian (de)stabilizations, see Figure \ref{lgmoves}. The proposition below follows the classical argument of \cite{LegendrianPaper}.

\bprop\label{isocusp} The links $\hat\kappa_1$ and $\hat\kappa_2$ are isotopic in $PT(S)$ if and only if $\kappa_1$ and $\kappa_2$ are related by Legendrian Reidemeister moves and Legendrian (de)stabilizations.
\eprop

\bpf The proof is identical to the $\bS^3$ case, see \cite[Theorem 4.4]{LegendrianPaper},
\epf

\begin{figure}[htb!]
\centering
\begin{overpic}[scale=0.5]{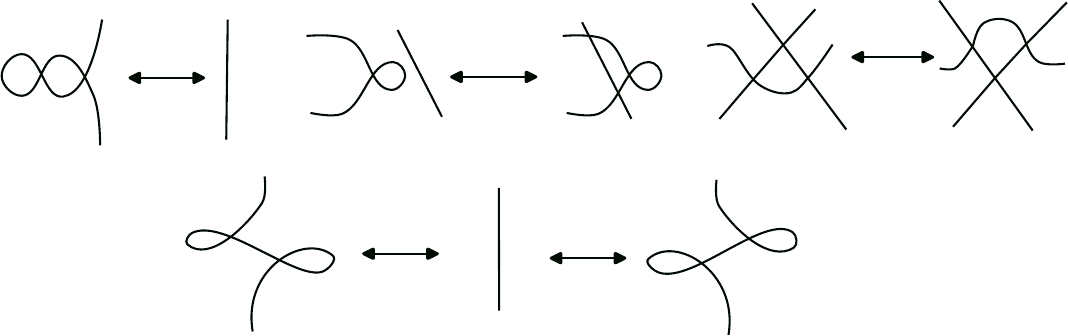}
\end{overpic}
\caption{Loop Reidemaister moves and (de)stabilizations. }\label{loopmoves}\end{figure}

Similarly, we can analyse the case for $UT(S)$. Here, loops essentially behave like cusps in this context. You can see the loop Reidemeister moves and loop (de)stabilizations in Figure \ref{loopmoves}. 

\bprop\label{isosmooth} The links $\hat\kappa_1$ and $\hat\kappa_2$ are isotopic in $UT(S)$ if and only if $\kappa_1$ and $\kappa_2$ are related by loop Reidemeister moves and loop (de)stabilizations.
\eprop

\bpf The proof is identical to the $\bS^3$ case, see \cite[Theorem 5.4]{LegendrianPaper},
\epf

Notice that transvections add cusps/loops while (de)stabilizations can both add or remove them. The cusps/loops added by transvections are harder to control. Theorem \ref{GL}, which we prove in Section \ref{canonproof}, states that the knot exterior of $\hat\kappa$ uniquely determines $\kappa$ up to cusps/loops added by transvections, Legendrian/loop Reidemeister moves and (de)stabilizations, and the action of mapping classes of $S$. In particular, if we restrict ourselves to the setting where $\kappa$ is smooth and in minimal position, then the equivalence is entirely given by mapping classes and Reidemeister III moves.

\subsection{Canonical links in tangent bundles}

To motivate our main application, let us begin by explaining how links in tangent bundles can be realized as canonical links. The idea is to generalize the classical fact about Lagrangian and Legendrian projections of links in $\mathbb{S}^3$.

\links*

\begin{proof} The proofs for $UT(S)$ and $PT(S)$ are identical, with the only change of loops replaces by cusps. Thus, we will focus on the $UT(S)$ variant. Fix an oriented link $\bar\kappa$ in $UT(S)$ an isotope it to be transverse to the fibers. Let $\kappa$ be the oriented shadow $\bar\kappa$ and apply a further isotopy so that $\kappa$ has transverse self-intersections and each intersection is between two arcs. Pick a complex structure and an abelian differential on $S$ such that $\kappa$ avoids the singularities of the differential. The abelian differential gives a horizontal and vertical foliation. Consider at the four induced vector fields, i.e. the left/right vector fields tangent to the horizontal and the up/down vector fields tangent to the vertical. We can trivialize $UT(S)$ away from the singularities and consider these vector fields as $\pi/2$-rotations of each other, with counter-clockwise rotation corresponding to going ``up'' in $UT(S)$.

We can now isotopy $\hat\kappa$ locally so that the shadow $\kappa$ is made up zig-zag paths switching between horizontal and vertical pieces with $\pi/2$ corners. Now, over each straight (horizontal or vertical) arc $\delta$, there is an annulus $A_\delta$ in $UT(S)$ and $\bar\kappa \cap A_\delta$ looks like $m$ positive or negative twists along this annulus. We can augment $\delta$ by adding $m$-loops, on the left or right, producing an arc $\delta_1$ with the property that the canonical lift $\hat\delta_1$ is isotopic to $\bar\kappa \cap A_\delta$ relative the points on $\partial A_\delta$. See figure \ref{addloops}. For crossings, we have a similar analysis, where we may have to add loops on either side so that the correct canonical arc is ``above'' or ``below,'' see Figure \ref{fixcrossing}. Lastly, we can smooth the corners by either smoothly turning or turning through a loop, as in Figure \ref{fixcorner}. The resulting augmented diagram on the surface if $\kappa_1$ and has the property that $\hat\kappa_1$ is isotopic to $\bar\kappa$ as desired.

\begin{figure}
\begin{center}
\begin{overpic}[scale=2]{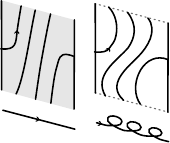}
\put(10,10){$\delta$}
\put(-10,70){$A_\delta$}
\put(28,80){$A_\delta \cap \bar\kappa$}
\put(60,5){$\delta_1$}
\put(60,60){$\hat\delta_1$}
\end{overpic}
\end{center}
\caption{On the left, we show $\delta$ and above $A_\delta$ with $\bar\kappa$. On the right, we show $\delta_1$ with its canonical lift $\hat\delta_1$. By convention, turning left goes ``up'' in $UT(S)$.}
\label{addloops}
\end{figure}

\begin{figure}
\begin{center}
\begin{overpic}[scale=2]{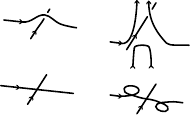}
\put(-20,15){$S$}
\put(-20,50){$UT(S)$}
\end{overpic}
\end{center}
\caption{On the left, we have a crossing shadow and above the associated piece of $\bar\kappa$. As drawn, the canonical lift of the shadow would have the opposite crossing. We correct is on the right by adding loops one either side. Note, the strands that go up glue to the strands entering form the bottom.}
\label{fixcrossing}
\end{figure}

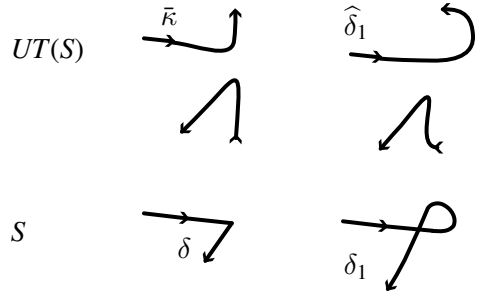
\begin{figure}
\begin{center}
\begin{overpic}[scale=2]{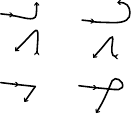}
\put(10,10){$\delta$}
\put(5,80){$\bar\kappa$}
\put(60,5){$\delta_1$}
\put(60,77){$\hat\delta_1$}
\put(-40,15){$S$}
\put(-40,70){$UT(S)$}
\end{overpic}
\end{center}
\caption{On the left, we have a corner shadow turning right, but $\bar\kappa$ goes up. On the right, we correct this with a loop. Note, the strands that go up glue to the strands entering form the bottom.}
\label{fixcorner}
\end{figure}

The case of $PT(S)$ is the same. The only necessary changes are to use a quadratic differential and to switch from loops to cusps. 
\end{proof}

We now see that this class of links very large. In fact, we have the following corollary.

\linkscor*

\bpf Recall that $M$ is an integral Dehn filling of $\mathbb{S}^3 \smallsetminus L$ for some link $L$ by \cite{L2012}. Thus, we can add  unknot $K$ to $L$ such that $L \subset \mathbb{S}^3 \setminus K$ can be seen as a link in a solid torus $N$. This solid torus $N$ can be realized as $UT(\mathbb{D}^2)$ or $PT(\mathbb{D}^2)$, allowing us to find a canonical diagram for $L$ in $N$. We now obtain $M$ as $1/0$-filling on $K$ and previous integral fillings on $L$.
\epf

\subsection{Exteriors of canonical links}\label{canonproof}

We now prove:

\appl*

\bpf Fist, we show the forward directions. Let $M = M_{\kappa_1}$. Since $M_{\kappa_1} \cong^+ M_{\kappa_2}$, we get two slopes $s_1, s_2$, corresponding to the meridians of $\hat{\kappa_1}, \hat{\kappa_2}$, such that the fillings of $M$ are the unit or projective bundles over $S_1, S_2$. By Lemma \ref{type}, $S_1 \cong^+ S_2$, so these bundles are diffeomorphic. Thus, $s_1$ and $s_2$ satisfy the conditions of Theorem \ref{SFone} and therefore $s_1 = s_2$. It follows that $\hat\kappa_1$ and $\hat\kappa_2$ are isotopic after applying a mapping class of the ambient circle bundle. By Section \ref{reiddiffcan} and Propositions \ref{isocusp} and \ref{isosmooth}, the result follows.

Similarly, the backwards direction follows from Propositions \ref{isocusp} and \ref{isosmooth} since the mapping class group of either tangent bundle takes meridians to meridians.\epf

\bibliographystyle{plain}
\bibliography{mybib}

\end{document}